\documentclass[12pt]{amsart}
\usepackage{amssymb,latexsym, amsmath, amsxtra}
\usepackage{graphicx}
\newtheorem{conjecture}{Conjecture}
\newtheorem{theorem}{Theorem}

\begin{document}
\title{Pair correlation and twin primes revisited}
 
\author{Brian Conrey}
\address{American Institute of Mathematics, 600 East Brokaw Rd., San Jose, CA 95112, USA and School of Mathematics, University of Bristol, Bristol BS8 1TW, UK}
\email{conrey@aimath.org}
\author{Jonathan P. Keating}
\address{School of Mathematics, University of Bristol, Bristol BS8 1TW, UK}
\email{j.p.keating@bristol.ac.uk}

\thanks{We gratefully acknowledge support under EPSRC Programme Grant EP/K034383/1
LMF: L-Functions and Modular Forms.  Research of the first author was also supported 
by the American Institute of Mathematics and by a grant from the National Science 
Foundation. JPK is grateful for the following additional support: a grant from the 
Leverhulme Trust, a Royal Society Wolfson Research 
Merit Award, and a Royal Society Leverhulme Senior Research Fellowship, and a grant from the Air Force Office of Scientific Research, Air Force Material Command, 
USAF (number FA8655-10-1-3088). He is also pleased to thank the American Institute of Mathematics for hospitality during a visit where this work started.}

\date{\today}

\begin{abstract} We establish a connection between the conjectural two-over-two ratios formula for the Riemann zeta-function and a conjecture concerning correlations of a certain arithmetic function.  Specifically, we prove that the ratios conjecture and the arithmetic correlations conjecture imply the same result.  This casts a new light on the underpinnings of the ratios conjecture, which previously had been motivated by analogy with formulae in random matrix theory and by a heuristic recipe.
\end{abstract}

\maketitle

Montgomery in his famous pair correlation paper [Mon] used heuristics based on 
the Hardy-Littlewood conjecture concerning the distribution of prime pairs [HL] to conclude that pairs of zeros of the Riemann zeta-function have the same scaled statistics, in the limit in which their height up the critical tends to infinity, as pairs of eigenvalues of 
large random Hermitian matrices (or of unitary matrices with Haar measure). Montgomery did not give the details of the calculation involving twin primes in his paper, but that calculation
has been repeated with variations several times in the literature, for example by Bolanz [Bol], Keating [K]
Goldston and Gonek [GG],  
and Bogomolny and Keating [BK1, BK2].  Goldston and Montgomery [GM] proved rigorously that the pair correlation conjecture is equivalent to an asymptotic formula for the variance of the number of primes in short intervals, and Montgomery and Soundararajan [MS] proved that this variance formula follows from the Hardy-Littlewood prime-pair conjecture, under certain assumptions.  

In a slightly different vein, Bogomolny and Keating [BK3, BK4] 
and later Conrey and Snaith [CSn] developed methods to give more precise estimates for the 
pair correlation (and higher correlations) of Riemann zeros. Bogomolny and Keating gave four different heuristic methods to accomplish this, while Conrey and Snaith used a uniform version of what is known as the ratios conjecture from which assumption they could rigorously derive this precise form of pair correlation.  All of these methods lead to the same formulae.
 
In this paper we reconsider this circle of ideas from yet another perspective, namely that of deriving
a form of the ratios conjecture from consideration of correlations between the values of a certain arithmetic function.  This provides a new perspective on the underpinnings of the ratios conjecture, which previously had been motivated by analogy with formulae in random matrix theory and by a heuristic recipe [CFKRS1, CFKRS2, CFZ].  This is similar to how, in a recent series of papers [CK1--4]  we have shown that moment conjectures previously developed using random matrix theory [KS, CFKRS2] may be recovered from correlations of divisor sums.

The twin prime conjectures are easily stated in terms of the von Mangoldt function $\Lambda(n)$ which is the generating function for $-\zeta'/\zeta$ (see, for example, [G]):
$$-\frac{\zeta'}{\zeta}(s)=\sum_{n=1}^\infty \frac{\Lambda(n)}{n^s}$$
or equivalently
\begin{eqnarray*}
\Lambda(n)=\left\{ \begin{array}{ll} \log p & \mbox{if $n=p^k$ for some prime $p$}\\
0& \mbox{otherwise}\end{array}\right.
\end{eqnarray*}
In the Conrey-Snaith approach zeros of $\zeta(s)$ are detected as poles 
of $\frac{\zeta'}{\zeta}(s)$ which in turn is realized via
$$\left. \frac{\zeta'}{\zeta}(s)= \frac{d}{d\alpha} \frac{\zeta(s+\alpha)}{\zeta(s+\gamma)}\right|_{\alpha=0\atop
\gamma=0}.$$
Passing to coefficients, we write
$$\mathcal I_{\alpha, \gamma}(s)=\sum_{n=1}^\infty \frac{I_{\alpha,\gamma}(n)}{n^s}=\frac{\zeta(s+\alpha)}{\zeta(s+\gamma)};$$
explicitly 
$$I_{\alpha,\gamma}(n)=\sum_{de = n} \frac{\mu(e)}{d^\alpha e^\gamma}.$$
Notice that
$$I_{\alpha,\gamma}(n)=n^rI_{\alpha+r,\gamma+r}(n)$$
for any $r$.
Also we have
$$\Lambda(n)=\left. - \frac{d}{d\alpha}I_{\alpha,\gamma}(n)\right|
_{\alpha=0\atop
\gamma=0}.$$

Here  we will investigate the averages 
$$\mathcal R_{\alpha,\beta,\gamma,\delta}(T):=\int_0^\infty \psi\left(\frac t T\right) \frac{\zeta(s+\alpha)\zeta(1-s+\beta)}{\zeta(s+\gamma)\zeta(1-s+\delta)}~dt$$
where $s=1/2+it$ and $\psi(z)$ is holomorphic in a strip around the real axis and decreases rapidly on the real axis.
Not surprisingly, $\mathcal R$  is related to averages of the (analytic continuation of the) Rankin-Selberg convolution
\begin{eqnarray*}
\mathcal{B}_{\alpha,\beta,\gamma,\delta}(s):=\sum_{n=1}^\infty \frac{I_{\alpha,\gamma}(n)
I_{\beta,\delta}(n)}{n^s}.
\end{eqnarray*}
In fact, 
the simplest case of the ratios conjecture asserts that   
\begin{eqnarray} \label{eqn:rat}
\mathcal R_{\alpha,\beta,\gamma,\delta}(T)=\int_0^\infty  \psi\left(\frac t T\right) \bigg(
\mathcal{B}_{\alpha,\beta\gamma,\delta}(1)+ \left(\frac{t}{2\pi}\right)^{-\alpha-\beta}
\mathcal{B}_{-\beta,-\alpha, \gamma,\delta}(1)
 \bigg) ~dt+O(T^{1-\eta})
 \end{eqnarray}
 for some $\eta>0$.   It is also not surprising that $\mathcal R$ is connected  to 
 weighted averages over $n$ and $h$ of
 $$  I_{\alpha,\gamma}(n)I_{\beta,\delta}(n+h).$$
It is this connection that we are elucidating.
Using the  the $\delta$-method it transpires that these weighted averages may 
be expressed in terms of 
 \begin{eqnarray*}&&
 \mathcal C_{\alpha,\beta,\gamma,\delta}(s):=\frac{1}{(2\pi i)^2}\int_{|w-1|=\epsilon}
 \int_{|z-1|=\epsilon}\chi(w+z-s-1)
 \sum_{q=1}^\infty 
 \sum_{h=1}^\infty \frac{r_q(h)}{h^{s+2-w-z}} \\
 &&\qquad  \qquad \qquad \qquad \times 
  \sum_{m=1}^\infty \frac{I_{\alpha,\gamma}(m)e(m/q)}{m^w}
 \sum_{n=1}^\infty \frac{I_{\gamma,\delta}(n)e(n/q)}{n^z}
 ~dw ~dz
  \end{eqnarray*}
  where $r_q(h)$ denotes Ramanujan's sum and where $\chi(s)$ is the factor 
  from the functional equation
  $\zeta(s)=\chi(s)\zeta(1-s)$; also here and elsewhere $\epsilon $ is chosen to be 
  larger than the absolute values of the shift parameters
  $\alpha,\beta,\gamma,\delta$ but smaller than $1/2$.
  The  result that ties this all together is the following identity.
    \begin{theorem}
  Assuming the Generalized Riemann Hypothesis
  \begin{eqnarray*}
  \mathcal C_{\alpha,\beta,\gamma,\delta}(s)=\mathcal B_{-\beta,-\alpha,\gamma,\delta}(s+1).
  \end{eqnarray*}
  \end{theorem}

In a recent series of papers [CK1--4] we have outlined a method that
 involves convolutions of coefficient correlations and leads to conclusions for averages 
 of truncations of products of shifted zeta-functions implied by the recipe
of [CFKRS2].  In this paper we strike out in a new direction, using similar ideas
to evaluate averages of truncations of products of ratios of shifted zeta-functions. 
In particular, the approach of Bogolmony \& Keating [BK1, BK2] on convolutions of shifted 
coefficient sums 
guide the calculations and we are led, as in the previous series, to formulate  
 a kind of  multi-dimensional Hardy-Littlewood circle method.
This first paper, as indicated above, may be viewed in a more classical context.

  It turns out to be convenient  to study an average of the ratios conjecture. To this end
let
$$\mathcal I_{\alpha, \gamma}(s;X)= \sum_{n\le X}
I_{\alpha,\gamma}(n)n^{-s}.
$$
We are interested in the average over $t$ of $\mathcal I_{\alpha,\gamma}\overline{\mathcal I_{\beta,\delta}}$  in the   case  that $X=T^\lambda$ for some $\lambda>1$ .
(When $\lambda<1$ this average is dominated by  diagonal terms.)
We give two different  treatments of the average of ``truncated'' ratios:
$$
\mathcal M_{\alpha,\beta,\gamma,\delta}(T;X):=\int_0^\infty \psi\left(\frac t T\right)  \mathcal I_{\alpha,\gamma}(s,X)
  \mathcal I_{\beta,\delta}(1-s,X)~dt$$
 (where again $s=1/2+it$) which lead to the same answer. The first is by the ratios conjecture and 
  the second is by consideration of the correlations of the coefficients.

  In each case we prove 
  \begin{theorem} Let $\alpha,\beta,\gamma,\delta$ be complex numbers smaller than $1/4$ 
  in absolute value.  Then, assuming either a uniform version of the ratios conjecture 
  or a uniform version of a conjectured formula for correlations of values of $I_{\alpha,\gamma}(n)$ (Conjecture 1, Section 4), we have for some $\eta>0$ and some $\lambda>1$,
  \begin{eqnarray*} &&
  \mathcal M_{\alpha,\beta,\gamma,\delta}(T;X)=\\
  &&   \qquad 
  \int_0^\infty  \psi\left(\frac t T\right) \frac{1}{2\pi i}\int_{\Re s=2}\bigg(
  \mathcal B_{\alpha,\beta,\gamma,\delta}(s+1)+\left(\frac t{2\pi}\right)^{-\alpha-\beta-s}
  \mathcal B_{-\beta,-\alpha,\gamma,\delta}(s+1)\bigg)\frac {X^s}{s}~ds~dt+O(T^{1-\eta}).
  \end{eqnarray*}
   \end{theorem}
  
  This shows that the ratios conjecture follows not only from the 'recipe' of [CFRKS2, CFZ], but also relates to  correlations of values of $I_{\alpha,\gamma}(n)$.
  
  \section{Approach via the ratios conjecture}
  
  We have
  $$\mathcal I_{\alpha,\gamma}(s,X)=\frac{1}{2\pi i} \int_{(2)} 
  \mathcal I_{\alpha,\gamma}(s+w)\frac{X^w}{w}~dw;
  $$
  there is a similar expression for $\mathcal I_{\beta,\delta}(s,X)$. 
  Inserting these expressions and rearranging the integrations we have 
  \begin{eqnarray*}
  \mathcal M_{\alpha,\beta,\gamma,\delta}(T;X)=\frac{1}{(2\pi i)^2}
  \int_{\Re w=2}\int_{\Re z=2} \frac{X^{w+z}}{wz} \mathcal R_{
\alpha+w,\beta+z,\gamma+w,\delta+z}(T)  ~dw ~dz.
\end{eqnarray*}
We observe from the expression (\ref{eqn:rat}) for the ratios conjecture that the integrand 
$\mathcal R_{\alpha+w,\beta+z,\gamma+w,\delta+z}$
 is, to leading order in $T$, expected to be a function of $z+w$.   
 We therefore make the change of variable $s=z+w$; now the integration in the $s$ variable is
 on the vertical line $\Re s=4$. We retain $z$ as our
 other variable and integrate over it. This turns out to be the integral
 $$\frac{1}{2\pi i} \int_{\Re z = 2}\frac{dz}{z(s-z)}=\frac{1}{s}$$
 as is seen by moving the path of integration to the left to $\Re z=-\infty$.
Thus we have that $\mathcal M_{\alpha,\beta,\gamma,\delta}(T;X)$ is given to leading order by
\begin{eqnarray*}
  \frac{1}{2\pi i}
  \int_{\Re s=4} \frac{X^{s}}{s} \mathcal R_{
\alpha+s,\beta,\gamma+s,\delta}(T)  ~ds.
\end{eqnarray*}
We move the path of integration to $\Re s=\epsilon$, avoiding crossing any poles, insert the ratios conjecture (\ref{eqn:rat}) (c.f.~the uniform version as laid out in 
[CSn]), and 
observe that  
$$\mathcal B_{\alpha+s,\beta,\gamma+s,\delta}(1)=\mathcal B_{\alpha,\beta,\gamma,\delta}(s+1).$$
In this way we have that the uniform ratios conjecture implies the conclusion of Theorem 2.

  \section{Approach via coefficient correlations}
  We follow the methodology developed in the work [GG] of Goldston and Gonek  
  on mean-values of long Dirichlet polynomials.

If we expand the sums and integrate term-by-term we have
$$\mathcal M_{\alpha,\beta,\gamma,\delta}(T;X)=T\sum_{m,n\le X} 
\frac{I_{\alpha,\gamma}(m)I_{\beta,\delta}(n)}{\sqrt{mn}}\hat\psi
\left(\frac{T}{2\pi}\log\frac mn\right).
$$
\subsection{Diagonal} The diagonal term is 
$$T\hat{\psi}(0) \sum_{m\le X} 
 \frac{I_{\alpha,\gamma}(m)I_{\beta,\delta}(m)}{m}.
 $$
 By Perron's formula the sum here is 
 \begin{eqnarray*}
 \frac{1}{2\pi i}\int_{(2)} \mathcal B_{\alpha,\beta,\gamma,\delta}(s+1)
 \frac{X^s}{s} ~ds.
 \end{eqnarray*}
 
\subsection{Off-diagonal}
 For the off-diagonal terms we  need to analyze
 $$2T\sum_{T\le m \le X}\sum_{1\le h \le \frac{X}{T}} \frac{I_{\alpha,\gamma}(m)
 I_{\beta,\delta}(m+h)}{m}\hat\psi
\left(\frac{Th}{2\pi m } \right).
$$
We replace the arithmetic terms by their average and express this as 
 $$2T\int_T^X \sum_{1\le h \le \frac{X}{T}} \frac{\langle I_{\alpha,\gamma}(m)I_{\beta,\delta}(m+h)\rangle _{m\sim u}
 }{u}\hat\psi
\left(\frac{Th}{2\pi u } \right) ~du.
$$
We compute the average  heuristically via the delta-method:
$$\langle I_{\alpha,\gamma}(m)I_{\beta,\delta}(m+h)\rangle _{m\sim u}
\sim \sum_{q=1}^\infty r_q(h) \langle I_{\alpha,\gamma}(m) e(m/q) \rangle_{m\sim u}
\langle I_{\beta,\delta}(m) e(m/q)  \rangle_{m\sim u}
$$
where $r_q(h)$ is the Ramanujan sum, a formula for which is $r_q(h)=
\sum_{d\mid h\atop d\mid q} d\mu(\frac q d )$; note that to actually prove this formula would be as difficult
as proving the Twin Prime conjecture.  We formalise this as a precise conjecture in Section 4.  It is this conjecture that we refer to in Theorem 2.
Now 
\begin{eqnarray*}
\langle I_{\alpha,\gamma}(m) e(m/q) \rangle_{m\sim u}
=\frac{1}{2\pi i}\int_{|w-1|=\epsilon}
\sum_{m=1}^\infty I_{\alpha,\gamma}(m)e(m/q) m^{-w} u^{w-1}~dw.
\end{eqnarray*}

Thus, the off-diagonal contribution is 
\begin{eqnarray*}&&
2T \sum_{1\le h \le \frac{X}{T}} \int_T^X
\frac{1}{(2\pi i)^2} \iint_{|w-1|=\epsilon\atop |z-1|=\epsilon} \sum_{q=1}^\infty r_q(h)
 \hat\psi
\left(\frac{Th}{2\pi u } \right)  u^{w+z-2} 
\\
&&\qquad \qquad \times \sum_{m_1=1}^\infty \frac{I_{\alpha,\gamma}(m_1)e(m_1/q)}{m_1^w}
\sum_{m_2=1}^\infty \frac{I_{\beta,\delta}(m_2)e(m_2/q)}{m_2^z}
~dw ~dz
~\frac{du}{u}.
\end{eqnarray*}
We make the change of variables $v=\frac{Th}{2\pi u }$. The inequality $u\le X$ then implies that 
$\frac{Th}{2\pi v} \le X$ or $h\le \frac{2\pi v X}{T}$.   The above can be re-expressed as
\begin{eqnarray*}&&
2T\int_0^\infty \sum_{1\le h \le \frac{2\pi vX}{T}} 
\frac{1}{(2\pi i)^2} \iint_{|w-1|=\epsilon\atop |z-1|=\epsilon} \sum_{q=1}^\infty r_q(h)
 \hat\psi(v) \left(\frac{Th}{2\pi v}\right)^{w+z-2} 
\\
&&\qquad \qquad \times \sum_{m_1=1}^\infty \frac{I_{\alpha,\gamma}(m_1)e(m_1/q)}{m_1^w}
\sum_{m_2=1}^\infty \frac{I_{\beta,\delta}(m_2)e(m_2/q)}{m_2^z}
~dw ~dz
~\frac{dv}{v}.
\end{eqnarray*}
Using Perron's formula to capture the sum over $h$ gives
\begin{eqnarray*}&&
2T\int_0^\infty  
\frac{1}{(2\pi i)^3}\int_{\Re s=2} \iint_{|w-1|=\epsilon\atop |z-1|=\epsilon} \sum_{q=1}^\infty 
\sum_{h=1}^\infty \frac{r_q(h)}{h^s}
 \hat\psi(v) \left(\frac{Th}{2\pi v}\right)^{w+z-2} \left(\frac {2\pi v X}T\right)^{s} 
\\
&&\qquad \qquad \times \sum_{m_1=1}^\infty \frac{I_{\alpha,\gamma}(m_1)e(m_1/q)}{m_1^w}
\sum_{m_2=1}^\infty \frac{I_{\beta,\delta}(m_2)e(m_2/q)}{m_2^z}\frac{ds}{s}
~dw ~dz
~\frac{dv}{v}.
\end{eqnarray*}
Now 
\begin{eqnarray*}
2\int_0^\infty \hat{\psi}(v) v^A \frac{dv}{v}=\chi(1-A)\int_0^\infty \psi(t) t^{-A} ~dt.
\end{eqnarray*}
Incorporating this formula leads us to 
\begin{eqnarray*}&&
T\int_0^\infty  \psi(t) 
\frac{1}{(2\pi i)^3}\int_{\Re s=2} \iint_{|w-1|=\epsilon\atop |z-1|=\epsilon} \sum_{q=1}^\infty 
\sum_{h=1}^\infty \frac{r_q(h)}{h^{s+2-w-z}}
  \left(\frac{Tt}{2\pi }\right)^{w+z-2} \left(\frac {2\pi X}{tT}\right)^{s} \chi(w+z-s-1)
\\
&&\qquad \qquad \times \sum_{m_1=1}^\infty \frac{I_{\alpha,\gamma}(m_1)e(m_1/q)}{m_1^w}
\sum_{m_2=1}^\infty \frac{I_{\beta,\delta}(m_2)e(m_2/q)}{m_2^z}\frac{ds}{s}
~dw ~dz
~dt.
\end{eqnarray*}
Hence, by Theorem 1, this is 
\begin{eqnarray*}&&
\int_0^\infty  \psi\left(\frac t T\right) 
\frac{1}{2\pi i}\int_{\Re s=2}  \left(\frac{t}{2\pi }\right)^{-\alpha-\beta-s}
 \mathcal B_{-\beta,-\alpha,\gamma,\delta}(s+1)\frac{X^s}{s}~ds
~dt.
\end{eqnarray*}

 Thus, adding the diagonal and off-diagonal terms 
 we obtain that the  conjecture for the correlations of values of $I_{\alpha,\gamma}(n)$ also implies the
 conclusion of Theorem 2.

\section{Proof of Theorem 1}

First of all, we have
$$\sum_{h=1}^\infty \frac{r_q(h)}{h^A} =\sum_{h=1}^\infty \frac{ \sum_{g\mid q\atop g\mid h}g
\mu(\frac q g)}{h^A}=\sum_{g\mid q} g^{1-A} \mu(\frac q g ) \zeta(A)
=q^{1-A}\Phi(1-A,q)\zeta(A)$$
where
$$\Phi(x,q)=\prod_{p\mid q}\left(1-\frac{1}{p^x}\right).$$
Using this and the functional equation for $\zeta$, we have to evaluate
\begin{eqnarray*}&&
\frac{1}{(2\pi i)^2} \iint_{|w-1|=\epsilon\atop |z-1|=\epsilon} \sum_{q=1}^\infty 
q^{w+z-s-1}\Phi(w+z-s-1,q)
   \\
&&\qquad \qquad \times
\zeta(w+z-s-1)
 \sum_{m_1=1}^\infty \frac{I_{\alpha,\gamma}(m_1)e(m_1/q)}{m_1^w}
\sum_{m_2=1}^\infty \frac{I_{\beta,\delta}(m_2)e(m_2/q)}{m_2^z} 
~dw ~dz.
\end{eqnarray*}

We can identify the polar structure of the Dirichlet series here by passing to characters
via the formula
$$e\left(\frac m q\right) = \sum_{d\mid m\atop d\mid q}
\frac{1}{\phi\left(\frac q d\right)}\sum_{\chi \bmod \frac q d } \tau(\overline{\chi})
\chi\left(\frac m d \right).$$
Assuming GRH, the only poles near $w=1$ arise from the principal characters 
$\chi_{\frac qd}^{(0)}$. Using 
$$\tau(\chi_{\frac qd}^{(0)})=\mu(\frac qd)$$
we have that  the poles of $\sum_{m=1}^\infty 
I_{\alpha,\gamma}(m)e(m/q) m^{-w}$ are the same as the poles of 
\begin{eqnarray*} &&
\sum_{d\mid q}\frac{\mu\left(\frac q d\right)}{\phi\left(\frac q d\right)}
 \sum_{m=1}^\infty I_{\alpha,\gamma}(md) \chi_{\frac qd}^{(0)}(m) m^{-w}d^{-w}\\
 &&\qquad = q^{-w}
 \sum_{d\mid q}\frac{\mu (d)}{\phi(d)}d^{w}
 \sum_{m=1}^\infty \frac{I_{\alpha,\gamma}(\frac{mq}d) \chi_{d}^{(0)}(m)}{ m^{w}} 
 \end{eqnarray*}
 and the principal parts are the same.
 We replace $\chi_d^{(0)}(m)$ by $\sum_{e\mid d\atop e\mid m}\mu(e)$.
 Thus we have
  \begin{eqnarray*}
 q^{-w}  \sum_{d\mid q}\frac{\mu (d)d^w}{\phi(d)}\sum_{e\mid d}\mu(e)e^{-w}
 \sum_{m=1}^\infty \frac{I_{\alpha,\gamma}(\frac{meq}d)  }{m^{w}} .
 \end{eqnarray*}
Now we need the polar structure of 
$$
 \sum_{m=1}^\infty I_{\alpha,\gamma}(mr)  m^{-w} 
 $$
 for $r=qe/d$.
 
 We use a lemma from [CGG] which asserts that if $A(w)=B(w)C(w)$ 
 where $A(w)=\sum_{m=1}^\infty \frac{a(m)}{m^{w}}$, $B(w)=\sum_{m=1}^\infty \frac{b(m)}{m^{w}}$
 and $C(w)=\sum_{m=1}^\infty \frac{c(m)}{m^{w}}$ then 
 \begin{eqnarray*}
 \sum_{m=1}^\infty \frac{a(mr)}{ m^{w}}=\sum_{r=r_1r_2}
 \sum_{m=1}^\infty \frac{b(mr_1)}{m^{w}}\sum_{m=1\atop (m,r_1)=1}^\infty
 \frac{c(mr_2)}{m^{w}}.
 \end{eqnarray*}
 We apply this identity with $a(m)=I_{\alpha,\gamma}(m)$, with $b(m)=m^{-\alpha}$ and 
 with $c(m)=\mu(m)m^{-\gamma}$ .
 Then
 $$\sum_{m=1}^\infty\frac{b(mr_1)}{m^w}=r_1^{-\alpha} \zeta(w+\alpha)$$
 and
  $$\sum_{(m,r_1)=1}\frac{c(mr_2)}{m^w}=
  \sum_{(m,r_1)=1}\frac{\mu(mr_2)}{m^{w+\gamma}r_2^\gamma}=
  \frac{\mu(r_2)}{r_2^{\gamma}}\sum_{(m,r)=1}\mu(m) m^{-w-\gamma}
  =\frac{\mu(r_2)r_2^{-\gamma}}{
  \Phi(w+\gamma,r)
   \zeta(w+\gamma)}.$$
   Now
    $$ \sum_{r=r_1r_2}\mu(r_2)r_1^{-\alpha}r_2^{-\gamma}
   =r^{-\alpha}\sum_{r=r_1r_2}\mu(r_2) r_2^{\alpha-\gamma}=r^{-\alpha}\Phi(\gamma-\alpha,r).$$
 Thus, 
  \begin{eqnarray*}
  \sum_{m=1}^\infty \frac{I_{\alpha,\gamma}(mr)}{  m^{w}}&=&
  \frac{\zeta(w+\alpha)  r^{-\alpha}\Phi(\gamma-\alpha,r)}{ \Phi(w+\gamma,r)\zeta(w+\gamma)} 
   \end{eqnarray*}
  
   In particular, we see that the only pole near to $w=1$ is at $w=1-\alpha$ with residue
    \begin{eqnarray*}
   \frac{ r^{-\alpha}
   \Phi(\gamma-\alpha,r)}{ \Phi(1+\gamma-\alpha,r)\zeta(1+\gamma-\alpha)} 
 . \end{eqnarray*}
 Inserting this with $r=qe/d$ into the above we now have that
 \begin{eqnarray*}
 \operatornamewithlimits{Res}_{w=1-\alpha}
 \sum_{m=1}^\infty\frac{I_{\alpha,\gamma}(m)e(\frac m q)}{m^w}&=& 
 q^{\alpha-1}  \sum_{d\mid q}\frac{\mu (d)d^{1-\alpha}}{\phi(d)}
 \sum_{e\mid d}\mu(e)e^{\alpha-1}
   \frac{ (qe/d)^{-\alpha} \Phi(\gamma-\alpha,qe/d)}{ \Phi(1+\gamma-\alpha,qe/d)
   \zeta(1+\gamma-\alpha)} \\
   &=& \frac{F_{\alpha,\gamma}(q)}{q\zeta(1+\gamma-\alpha)}
 \end{eqnarray*}
 where 
 $$F_{\alpha,\gamma}(q)=
 q^{\alpha}  \sum_{d\mid q}\frac{\mu (d)d^{1-\alpha}}{\phi(d)}\sum_{e\mid d}\mu(e)e^{\alpha-1}
   \frac{ (qe/d)^{-\alpha} \Phi(\gamma-\alpha,qe/d)}{ \Phi(1+\gamma-\alpha,qe/d) } $$ 
 is a multiplicative function of $q$.
 At a prime $p$ we have
 \begin{eqnarray*}
 F_{\alpha,\gamma}(p)&=&
 p^{\alpha} \bigg(\frac{p^{-\alpha}\Phi(\gamma-\alpha,p)}{\Phi(1+\gamma-\alpha,p)}
 -\frac{p^{1-\alpha}}{p-1}\big(1-\frac{p^{\alpha-1}p^{-\alpha}\Phi(\gamma-\alpha,p)}
 {\Phi(1+\gamma-\alpha,p)}\big)\bigg)\\
 &=&\frac{ \Phi(\gamma-\alpha,p)}{\Phi(1+\gamma-\alpha,p)}\left(1+\frac{1}{p-1}\right)
 -\frac{p}{p-1}\\
 &=& \frac{p}{(p-1)}\left(\frac{ \Phi(\gamma-\alpha,p)}{\Phi(1+\gamma-\alpha,p)}-1\right)
 =\frac{p}{(p-1)}\left(\frac{ (1-p^{\alpha-\gamma})}{
 (1-p^{-1+\alpha-\gamma})}
 -1\right)\\
 &=& \frac{p}{(p-1)}  \frac{ (-p^{\alpha-\gamma}+p^{-1+\alpha-\gamma})}{
 (1-p^{-1+\alpha-\gamma})}= \frac{- p^{\alpha-\gamma} }{
 (1-p^{-1+\alpha-\gamma})}
 =-p^{\alpha-\gamma}+O(\frac 1 p ).
 \end{eqnarray*}

 With $w=1-\alpha$ and $z=1-\beta$ we see that our sum is 
 \begin{eqnarray*}&&
\frac{\zeta(1-\alpha-\beta-s)}{\zeta(1-\alpha+\gamma)\zeta(1-\beta+\delta)} 
 \sum_{q=1}^\infty 
q^{-1-\alpha-\beta-s}\Phi(1-\alpha-\beta-s,q) F_{\alpha,\gamma}(q) F_{\beta,\delta}(q)
   \end{eqnarray*}
Because of  $  F_{\alpha,\gamma}(p) = 
  -p^{\alpha-\gamma}+O(\frac 1 p )$ we have
  \begin{eqnarray*}
\sum_{q=1}^\infty 
q^{-1-\alpha-\beta-s}\Phi(1-\alpha-\beta-s,q) F_{\alpha,\gamma}(q) F_{\beta,\delta}(q)
=\zeta(1+\gamma+\delta+s)B_{\alpha,\beta,\gamma,\delta}(s)
\end{eqnarray*}
where $B$ is an Euler product that is absolutely convergent for $s$ near 0. 
We claim that $B_{\alpha,\beta,\gamma,\delta}(s)=A_{-\beta, -\alpha-s,\gamma+s,\delta}.$
 This is easily seen to be equivalent to showing that 
$$ B_{\alpha,\beta,\gamma,\delta}(0)=A_{-\beta, -\alpha,\gamma,\delta}.$$
To prove this we first note that for $j\ge 2$ we have 
 \begin{eqnarray*}
F_{\alpha,\gamma}(p^j)&=&p^{j\alpha}\bigg(\frac{p^{-j\alpha}\Phi(\gamma-\alpha,p)}
{\Phi(1+\gamma-\alpha,p)}
-\frac{p^{1-\alpha}}{p-1}
\bigg(\frac{p^{-(j-1)\alpha}
\Phi(\gamma-\alpha,p)}{\Phi(1+\gamma-\alpha,p)}-p^{\alpha-1}
\frac{p^{-\alpha j}\Phi(\gamma-\alpha,p)}{\Phi(1+\gamma-\alpha,p)}
\bigg)\bigg)
\\
&=&\frac{\Phi(\gamma-\alpha,p)}{\Phi(1+\gamma-\alpha,p)}\bigg(
1-\frac p{(p-1)} +p^{\alpha-1}
\bigg)=\frac{\Phi(\gamma-\alpha,p)}{\Phi(1+\gamma-\alpha,p)}\bigg(
 -\frac 1{(p-1)} +\frac 1{(p-1)}
\bigg)=0.
\end{eqnarray*}
 
 Now the sum of the series 
 $$ \sum_{j=0}^\infty 
p^{(-1-\alpha-\beta)j}\Phi(1-\alpha-\beta,p^j) F_{\alpha,\gamma}(p^j) F_{\beta,\delta}(p^j)
$$ is just
\begin{eqnarray*}&&
 1+p^{-1-\alpha-\beta}\Phi(1-\alpha-\beta,p) F_{\alpha,\gamma}(p) F_{\beta,\delta}(p)\\&&\qquad = 
1+\frac{(1-\frac{1}{p^{1-\alpha-\beta}})}{p^{1+\alpha+\beta}}\frac{ p^{\alpha-\gamma} }{
 (1-p^{-1+\alpha-\gamma})}\frac{ p^{\beta-\delta} }{
 (1-p^{-1+\beta-\delta})}
 \\&&\qquad=
 1+\frac{(1-\frac{1}{p^{1-\alpha-\beta}})}
{ p^{1+ \gamma+\delta}
 (1-p^{-1+\alpha-\gamma})(1-p^{-1+\beta-\delta})}\\
 &&\qquad = (1-\frac{1}{p^{1+\gamma+\delta}})^{-1}B^{(p)}_{\alpha,\beta,\gamma,\delta}(0)
\end{eqnarray*}
where
\begin{eqnarray*}
  B^{(p)}_{\alpha,\beta,\gamma,\delta}(0)&=&
   (1-\frac{1}{p^{1+\gamma+\delta}})\left(
    1+\frac{(1-\frac{1}{p^{1-\alpha-\beta}})}
{ p^{1+ \gamma+\delta}
 (1-p^{-1+\alpha-\gamma})(1-p^{-1+\beta-\delta})}\right)
  \end{eqnarray*}
  The identity will be proven provided we can show that
  \begin{eqnarray*}
   1+\frac{(1-\frac{1}{p^{1-\alpha-\beta}})}
{ p^{1+ \gamma+\delta} 
 (1-p^{-1+\alpha-\gamma})(1-p^{-1+\beta-\delta})}= \frac{(1-\frac{1}{p^{1-\alpha+\gamma}}
 -\frac{1}{p^{1-\beta+\delta}}+\frac{1}{p^{1+\gamma+\delta}})}
 {(1-\frac{1}{p^{1-\beta+\delta}})(1-\frac{1}{p^{1-\alpha+\gamma}})}
 \end{eqnarray*}
 This is equivalent to showing that
 \begin{eqnarray*}
 1+\frac{XCD(1-\frac{X}{AB})}{(1-\frac{XC}{A})(1-\frac{XD}{B})}
 = \frac{(1-\frac{XC}{A}-\frac{XD}{B}+XCD)}{(1-\frac{XD}{B})(1-\frac{XC}{A})}
 \end{eqnarray*}
 where $X=\frac{1}{p}$; $A=p^{-\alpha}$; $B=p^{-\beta}$; $C=p^{-\gamma}$;  $D=p^{-\delta}$.
 This reduces to 
 \begin{eqnarray*}
 (1-\frac{XC}{A})(1-\frac{XD}{B})+XCD(1-\frac{X}{AB})
 =(1-\frac{XC}{A}-\frac{XD}{B}+XCD)
 \end{eqnarray*}
 or
  \begin{eqnarray*}
 (A-XC)(B-XD)+XCD(AB-X ) 
 =AB -XC -XD +XABCD 
 \end{eqnarray*}
 which is easily checked.
 
 \section{ Conjecture 1}
 
 We can use the results of the previous two sections to formulate the  conjecture that is part of the input 
 for Theorem 2.
 
 We expect $I_{\alpha,\gamma}(n)I_{\beta,\delta}(n+h)$ 
 for $n$ near $u$ to behave on average like
 $$\sum_{q=1}^\infty r_q(h) \frac{1}{(2\pi i)^2}
 \int_{|w-1|=\epsilon}
 \sum_{m=1}^\infty \frac {I_{\alpha,\gamma}(m)e(m/q)}{m^w} u^{w-1} ~dw
 \int_{|z-1|=\epsilon} 
  \sum_{n=1}^\infty \frac {I_{\beta,\delta}(n)e(n/q)}{n^z} u^{z-1} ~dz.
 $$
 The integrals over $w$ and $z$ are 
 $$\frac{F_{\alpha,\gamma}(q)u^{-\alpha}}{q\zeta(1+\gamma-\alpha)}
 \qquad
 \frac{F_{\beta,\delta}(q)u^{-\beta}}{q\zeta(1+\delta-\beta)}$$
 respectively. Thus, $I_{\alpha,\gamma}(n)I_{\beta,\delta}(n+h)$ 
behaves like 
$$  \frac{n^{-\alpha-\beta}}{\zeta(1+\gamma-\alpha)\zeta(1+\delta-\beta)}
     \sum_{q=1}^\infty \frac{r_q(h)  F_{\alpha,\gamma}(q) F_{\beta,\delta}(q) }{q^2}
 .
 $$
In particular, we expect that
$$\sum_{n=1}^\infty \frac{I_{\alpha,\gamma}(n)I_{\beta,\delta}(n+h)}{n^s}-
 \frac{\zeta(s+\alpha+\beta)}{\zeta(1+\gamma-\alpha)\zeta(1+\delta-\beta)}
     \sum_{q=1}^\infty \frac{r_q(h)  F_{\alpha,\gamma}(q) F_{\beta,\delta}(q) }{q^2}
$$
is analytic in $\sigma>\sigma_0$ for some $\sigma_0< 1$.

This leads us to 
\begin{conjecture}There are numbers $\phi<1$
and $\psi>0$ such that 
\begin{eqnarray*}
\sum_{n\le x}I_{\alpha,\gamma}(n)I_{\beta,\delta}(n+h)=m(x,h)+O(x^\phi)
\end{eqnarray*}
uniformly for $h\ll x^\psi$ where 
$$m(x,h)= \frac{1}{\zeta(1+\gamma-\alpha)\zeta(1+\delta-\beta)}
     \sum_{q=1}^\infty \frac{r_q(h)  F_{\alpha,\gamma}(q) F_{\beta,\delta}(q) }{q^2}\frac{x^{1-\alpha-\beta}}{1-\alpha-\beta}.$$
\end{conjecture}
 \section{Conclusion}
 
 In subsequent papers we will extend this process
 to averages of truncated ratios with any number of factors in the numerator and denominator.

 \section{Appendix}
 For ease of comparison with results in the literature we 
 give a more concrete expression for $\mathcal M$.
 
  First of all, we note that 
  the Rankin-Selberg Dirichlet series has an Euler product
 \begin{eqnarray*}
 \mathcal B_{\alpha,\beta,\gamma,\delta}(s)=\sum_{m=1}^\infty  \frac{I_{\alpha,\gamma}(m)I_{\beta,\delta}(m)}{m^{s}}
 =\prod_p \sum_{j=0}^\infty  \frac{I_{\alpha,\gamma}(p^j)I_{\beta,\delta}(p^j)}{p^{js}}.
 \end{eqnarray*}
 Now
 $$\sum_{j=0}^\infty I_{\alpha,\gamma}(p^j) x^j =\frac{1-p^{-\gamma}x}{1-p^{-\alpha}x}
 =(1-p^{-\gamma}x)(1+p^{-\alpha}x+p^{-2\alpha}x^2 + \dots)$$
 so that 
 \begin{eqnarray*}I_{\alpha,\gamma}(p^j)=\left\{ \begin{array}{ll} p^{-\alpha j} (1-p^{\alpha-\gamma})& \mbox{if $j\ge 1$}\\
 1 &\mbox{if $j=0$}\end{array} \right.  
 \end{eqnarray*}
 Thus,
 \begin{eqnarray*}
 \sum_{j=0}^\infty I_{\alpha,\gamma}(p^j) I_{\beta,\delta}(p^j) x^j&=&
 1+ (1-p^{\alpha-\gamma})(1-p^{\beta-\delta})\sum_{j=1}^\infty p^{-(\alpha+\beta) j} x^j\\
 &=& \frac{1- p^{-\beta-\gamma}x -p^{-\alpha-\delta}x +p^{-\gamma-\delta}x}
 {1-p^{-\alpha-\beta}x}
 \end{eqnarray*}
 and 
  \begin{eqnarray*}
 \sum_{m=1}^\infty  \frac{I_{\alpha,\gamma}(m)I_{\beta,\delta}(m)}{m^{s}}
 &=&\zeta(s+\alpha+\beta)\prod_p
\left(1- \frac{1}{p^{s+\beta+\gamma}}
 -\frac{1}{p^{s+\alpha+\delta}} +\frac{1}{p^{s+\gamma+\delta}}\right)\\
 &=& \frac{\zeta(s+\alpha+\beta)\zeta(s+\gamma+\delta)}
 {\zeta(s+\alpha+\delta)\zeta(s+\beta+\gamma)}A_{\alpha,\beta,\gamma,\delta}(s)
 \end{eqnarray*} 
 where  \begin{eqnarray*}
 A_{\alpha,\beta,\gamma,\delta}(s)=\prod_p\frac{
 \left(1-\frac{1}{p^{s+\gamma+\delta}}\right)\left(1- \frac{1}{p^{s+\beta+\gamma}}
 -\frac{1}{p^{s+\alpha+\delta}} +\frac{1}{p^{s+\gamma+\delta}}\right)}
 {\left(1-\frac{1}{p^{s+\beta+\gamma}}\right)\left(1-\frac{1}{p^{s+\alpha+\delta}}\right)}.
 \end{eqnarray*}

Now it is an easy exercise to calculate that
 \begin{eqnarray*} &&
  \mathcal M_{\alpha,\beta,\gamma,\delta}(T;X)=\\
  &&   \qquad 
  \int_0^\infty  \psi\left(\frac t T\right) \bigg(
\frac{\zeta(1+\alpha+\beta)\zeta(1+\gamma+\delta)}
 {\zeta(1+\alpha+\delta)\zeta(1+\beta+\gamma)}A_{\alpha,\beta,\gamma,\delta}(1)\\
 &&\qquad \qquad + \left(\frac{t}{2\pi}\right)^{-\alpha-\beta}
 \frac{\zeta(1-\beta-\alpha)\zeta(1+\gamma+\delta)}
 {\zeta(1-\beta+\delta)\zeta(1-\alpha+\gamma)}A_{-\beta,-\alpha ,\gamma,\delta}(1)
 \\&&\qquad \qquad -
 \frac{X^{-\gamma-\delta}}{(\gamma+\delta)}
\frac{ \zeta(1+\alpha+\beta-\gamma-\delta)}
 {\zeta(1 +\alpha -\gamma )\zeta(1+\beta -\delta) }A_{\alpha-\gamma-\delta,\beta,-\delta,
 \delta}(1)\\
 &&\qquad  \qquad +  \left(\frac{t}{2\pi}\right)^{-\alpha-\beta}
 \left(\frac{t}{2\pi X}\right) ^{\gamma+\delta} \frac{ \zeta(1+\gamma+\delta-\alpha-\beta)}
 {\zeta(1-\alpha+\gamma) \zeta(1-\beta+\delta) (\gamma+\delta)}
A_{-\beta,\gamma+\delta-\alpha,-\delta,\delta}(1)
 \bigg) ~dt\\ &&\qquad \qquad \qquad +O(T^{1-\eta})
 \end{eqnarray*}
for some $\eta >0$.

\end{document}